\documentclass[authoryear,12pt,3p]{jowarticle}

\usepackage[english]{babel}
\usepackage[utf8]{inputenc}
\usepackage{amsmath}
\usepackage{graphicx}
\usepackage{natbib}
\usepackage{blindtext}
\usepackage{subcaption}
\usepackage{setspace}

\makeatletter
\renewcommand\section{\@startsection{section}{1}{\z@}{-3.25ex plus -1ex minus -.2ex}{1.5ex plus .2ex}{\normalsize\bf}}
\renewcommand\subsection{\@startsection{subsection}{2}{\z@}{-3.25ex plus -1ex minus -.2ex}{1.5ex plus .2ex}{\normalsize\bf}}
\renewcommand\subsubsection{\@startsection{subsubsection}{3}{\z@}{-3.25ex plus -1ex minus -.2ex}{1.5ex plus .2ex}{\normalsize\bf}}
\makeatother

\begin{document}
\begin{frontmatter}
\title{Correctness, Artificial Intelligence, and the\\ Epistemic Value of Mathematical Proof}
\author{James Owen Weatherall}\ead{james.owen.weatherall@uci.edu} 
\address{Department of Logic and Philosophy of Science \\ University of California, Irvine}
\author{Jesse Wolfson}\ead{wolfson@uci.edu} 
\address{Department of Mathematics \\ University of California, Irvine}

\date{\today}

\begin{abstract}
We argue that it is neither necessary nor sufficient for a mathematical proof to have epistemic value that it be ``correct'', in the sense of formalizable in a formal proof system.  We then present a view on the relationship between mathematics and logic that clarifies the role of formal correctness in mathematics.  Finally, we discuss the significance of these arguments for recent discussions about automated theorem provers and applications of AI to mathematics.
\end{abstract}
\end{frontmatter}

\doublespacing
\section{Introduction} \label{sec:intro}

That a mathematical proof be correct is neither sufficient nor necessary for it to have epistemic value.  By ``correct'', here, for the sake of the argument, we mean that the proof can be formalized, i.e., translated into a formal proof system.  (Below, we will call this sense of correctness ``formal correctness'' and distinguish it from ``mathematical correctness''.)  By having ``epistemic value'', we mean being well-suited to play the various epistemic roles that proof does, and should, play in mathematical practice.\footnote{Some readers will think that by ``epistemic value'' we are primarily interested in \emph{justification}, and that our arguments are aimed at showing that something proof-like (a ``simil-proof'' in the terminology of \citet{deToffoliGFM}, which is to say, apparent proofs lacking obvious errors that may nonetheless later be shown to be fallacious) can provide justification, even if it is not formally correct \citep[c.f.][]{DeToffoliJWP}.  But as we hope will be clear in what follows, we intentionally do not center our discussion around justification, and we take the epistemic value of proofs to be more general.} In other words, epistemic value is what makes proofs ``work'', as proofs, for the purposes that mathematicians seek to attain.  We will say more below about what these roles are, and about how we see the relationship between formal and informal proofs.  

There are two related strands of literature that form the context for the present arguments.  One is the philosophical literature on the Standard View of mathematical proof and rigor \citep{Azzouni, Avigad2006, Avigad2008, Avigad2021, Avigad2022, Weir, HamamiInference, HamamiRigour, Tatton-Brown}, which holds that ordinary ``informal'' proofs, as generated by working mathematicians and published in high quality mathematics journals, should be seen as describing, or encoding, or ``indicating'' \citep{Azzouni}, or ``sketching'' \citep[p. 377]{Maclane} a formal proof that could, at least in principle, be generated from the informal proof (though whether this process is ``routine'' \citep{Maclane} or an additional ``creative act'' \citep{Avigad2021,Larvor} varies among proponents of the view).  The other strand of literature -- or at least, conversation -- has appeared mostly in talks and online forums,\footnote{We have in mind examples like posts on Kevin Buzzard's blog Xena, Michael Harris's blog Mathematics without Apologies, frequent and influential comments by Terence Tao on Mastodon, various guest blog posts from prominent mathematicians \citep[e.g.][]{Gowers,Tao,Scholze}, and widely viewed online talks by, for instance, \citet{Buzzard2022,Buzzard2023} and \citet{TaoIMO,TaoSimons}, among others.  See also the interview Tao gave to \emph{Scientific American} on the future of AI in mathematics \citep{SciAm}.} with some cross-over to mathematics journals \citep[e.g][]{Fraser+etal,Venkatesh}, concerning the speculative possibility of AI systems that could generate correct proofs of theorems beyond the ken of human mathematicians, much like AI systems have exceeded human abilities in games like Go and Chess. 

Our thesis may give readers pause.  Of course mathematical proofs should be correct!  And if a proof cannot be formalized, surely that is a signal that something odd or surprising is going on.  We grant this---as, we think, would most contributors to the literature on the Standard View.  But even so, there is a long history of debate over the status of the Standard View, arguably originating in debates between figures such as Hilbert, Brouwer, and Poincar\'e in the early 20th century, with recent entries in the philosophical literature from \citet{Rav1999}, \citet{Detlefsen}, \citet{Cellucci}, \citet{LarvorOld}, and \citet{Tanswell} criticizing aspects of the Standard View.\footnote{A rich and parallel discourse in the mathematics community runs through \citet{Polya}, \citet{Lakatos}, \citet{Thurston}, \citet{HershBook}, \cite{Manin}, \citet{Arnold} and \citet{HarrisQuanta,HarrisDogma}.}  (We will discuss some of these other arguments in more detail below.)  These arguments tend to have one of two forms: either they emphasize that formal proofs are purely syntactic, whereas informal ones seem to have some further content related to knowledge of the subject matter, so that they are not \emph{merely} abbreviations of formal proofs; or else they emphasize that there are significant mismatches between what appears to be required for a successful formal proof and what is needed for an informal one.

The present paper should be seen as part of the same tradition.  But we wish to put our emphasis in a different place---and to offer what we believe is a deeper critique of the Standard View, grounded in a form of anti-logicism that we think is prevalent in many areas of mathematics, but which is not often expressed or defended in the literature on the Standard View.\footnote{That said, we do not claim the view is original in philosophy of mathematics---see, for instance, \citet{MaddyEIT}, or related ideas in \citet{BurgessProofs}.}  On our view, while it seems to be generally true that mathematically correct informal proofs can be made formally correct, this fact is incidental and secondary to the epistemic value of informal proof.  Proofs are formally correct, when they are, not because correctness is the normative ideal governing mathematicians' activities, nor because correctness underwrites the epistemic value of those activities.  Instead, informal proofs are formally correct because logicians have developed a rich and successful mathematical theory of mathematical practice intended to capture various aspects of mathematical argument.  Formal correctness tracks mathematical correctness because formal correctness was designed to do precisely that.  

In other words, the formal systems in which formal correctness is established have been invented and tuned with the goal of making existing informal proof methods formalizable.  Formalizability speaks to the success of this activity in mathematical logic---but it is irrelevant to the question of what makes the informal proofs epistemically valuable.  Insofar as failures of correctness signal problems with proofs, it is not because formal correctness is the goal of the proof; rather, it is that failures of formal correctness indicate salient differences in the reasoning between the (formally) incorrect proof and other examples of successful mathematical reasoning.  Conversely, the correctness of a proof is not what gives it epistemic value, even \emph{qua} proof. 

Although our motivations for criticizing the Standard View are similar to those of previous authors, we feel the issue has become increasingly urgent.  In recent years, we have seen rapid advances in the areas of computer-assisted theorem proving and proof verification.  Those working in these areas often invoke the ideology of the Standard View to explain and justify their goals: if it is correctness that is the ultimate goal, then automated systems that can ensure correctness are the ideal proof generators.  This sort of reasoning has led to some striking soothsaying: in the future, one should expect AI mathematicians to outpace human mathematicians at proving theorems, in much the same way that AI game players have outpaced their human competitors in Go, chess, Jeopardy, and other games \citep{Bory}. Many groups have already devoted substantial resources to realizing this vision of future mathematics; and some very prominent mathematicians have become boosters for AI-driven mathematics (e.g. \cite{Gowers,Tao}).

We do not want to gainsay the future value of AI for mathematics.  We, too, are impressed by the rate of progress in recent years, and we suspect that computer assistance will change many aspects of mathematical life---just as Mathematica and Matlab changed the lives of undergraduate calculus students and many applied mathematicians a generation ago.  But we maintain that if one accepts our principal theses, then the mere generation of correct proofs is not, in itself, a contribution to mathematics.\footnote{To be sure, many advocates of the Standard View can and will agree with this observation.  It is nonetheless important to repeat, especially in the context of understanding what AI mathematics can contribute.}  Something more is needed. The broader mathematical community needs to be clear about what these sorts of activities can and cannot contribute, and what risks are posed to more traditional mathematical activities by reorienting towards machine-driven methods.  While we will not attempt to give an exhaustive list of what additional criteria are jointly sufficient for proofs to have epistemic value, we will point below to several conditions that we think are necessary.  

The remainder of the paper will proceed as follows.  We will begin by reviewing the Standard View and some previous criticisms of it to better situate what we are trying to do here.  In the following section, we will elaborate on what we mean by the ``epistemic value'' of a proof, and use that machinery to argue that correctness cannot be sufficient for a proof to have epistemic value.  We will then turn to the other half of our claim, and argue that correctness is not necessary.  Here we will offer several examples of cases in which incorrect attempts at proofs (or incorrect simil-proofs) turned out to be of exceptional epistemic value.  Then, in section \ref{sec:anti-logicism} we will return to the view sketched above of how logic should be seen as related to mathematics, and we will offer an account of why one might have thought formal correctness was necessary.  We will conclude by reflecting on how our arguments bear on the question that motivated us, concerning what sort of contributions AI systems can make to the mathematical enterprise.

\section{The Standard View: For and Against}

\subsection{For the Standard View}

The Standard View of mathematical proof is a claim about the relationship between mathematical proof as encountered in the wilds of mathematical practice and formal proofs.\footnote{The present review is intended only to give a sense of the literature and motivate the arguments that follow. For an extended introduction and exhaustive literature review, see \citet{TanswellElement}; see also \citet{Burgess+DeToffoli}.}  One of the Standard View's most prominent defenders describes it as follows:
\begin{quote} \singlespacing When someone in the mathematical community makes a mathematical claim, it is generally possible to express that claim formally, in the sense that logically adept and sufficiently motivated mathematicians can come to agreement that the formal claim expresses the relevant theorem.  One justifies an informal claim by proving it, and if the proof is correct, with enough work it can be turned into a formal derivation.  Conversely, a formal derivation suffices to justify the informal claim.  So an informal mathematical statement is a theorem if and only if its formal counterpart has a formal derivation.  Whether or not a mathematician reading a proof would characterize the state of affairs in these terms, a judgement as to correctness is tantamount to a judgment as to the existence of a formal derivation, and whatever psychological processes the mathematician brings to bear, they are reliable insofar as they track the correspondence. \citep[p. 7379]{Avigad2021}\end{quote}
In other words, a proof plays a justificatory role in establishing a theorem only if the proof can be translated into a formal derivation; and insofar as a proof is ``correct'', it \textit{can} be so translated.\footnote{Per personal communication, Avigad himself would prefer a weaker term than ``translated''. (Still, we think that ``translated'' is a reasonable way to summarize ``with enough work can be turned into''.)  His view is that actually generating a formal proof may involve filling gaps, correcting errors, or other work that may extend beyond ``translation''.}

The Standard View as we have just described it is hardly new.  Some authors identify its source in Frege, though it is not clear if he truly held a version of the Standard View, since he maintained every inference has some ``non-formal'' component.  But there are good reasons to think that both Hilbert and members of his school subscribed to the Standard View.  \citet{Cellucci}, for instance, calls the view the Hilbert-Gentzen Thesis \citep[c.f.][]{Hilbert,Gentzen}, and states it as follows: ``Every real proof can be represented by a formal proof.''  \citet{HamamiRigour} highlights articulations of the view by \citet{Maclane} and \citet{Bourbaki}, though he also argues that it is broadly adopted by 20th century mathematics after Hilbert.   It is also endorsed by some contemporary mathematicians.  And as we have noted, versions of the Standard View appear to be behind some mathematicians'  calls for a significant investment in automated proof checking or  proof generation (e.g. \cite{Hales2008,Grayson2017,Scholze,Buzzard2024}).

Careful modern philosophical articulations of the Standard View come in several flavors.\footnote{There are also variations in whether the Standard View concerns proof, \emph{per se}, or merely \emph{rigor} as a properties of proofs.}  On some versions of the view, an informal proof contains the information needed to construct the formal one.  For instance, \citet[p. 409]{HamamiRigour} argues that ``Mathematical proof is the primary form of justification of mathematical knowledge.  But in order to count as a \textit{proper} mathematical proof, and thereby to function \textit{properly} as a justification for a piece of mathematical knowledge, a mathematical proof must be \textit{rigorous}.''  He goes on to describe (and defend) what he calls the Standard View: ``According to this view, a mathematical proof $P$ is \textit{rigorous} if and only if $P$ can be \textit{routinely translated} into a formal proof (p. 410).''   

Other versions are more abstract.  \citet{Azzouni}, for instance, defends what he calls the ``derivation-indicator view of ordinary mathematical proof'' (p. 85), which is the view that the role of proof in practice is to convince other mathematicians of the existence of a certain kind of ``derivation'', or formal proof---and not, necessarily, to directly encode that proof.  This is a version of the standard view that elaborates the reasons -- pragmatics, efficiency, expedience -- that mathematicians give informal proofs, even though it is ultimately formal proofs that matter.  For Azzouni, the Standard View has important explanatory virtues: in particular, it explains why there is intersubjective agreement among mathematicians about which proofs are correct.  It is because these proofs are the ones that successfully indicate the existence of a formal proof.

\subsection{Intermezzo}

Before proceeding to review arguments against the Standard View, let us pause to clarify what we intend to accomplish in what follows in the light of the previous subsection.  For our purposes, the differences between the various versions of the Standard View will not matter.  We see ``(formally) correct'', ``rigorous'', ``proper'', and other similar terms as success terms indicating that a proof has achieved its key epistemic aim -- usually, in this literature, providing justification for a proposition -- and therefore has epistemic value in the sense described above, i.e., that it is able to play the epistemic role in mathematical practice that proof is expected to play.  And so we take the Standard View to be the claim that the epistemic value of an informal proof follows from the in-principle existence of a correct formal proof that bears some relationship (``formalization'') to it.  Our arguments are intended to respond to this broad family of views.

Some readers will see a gap, or perhaps a change of subject, between the Standard View and the main claims we articulated in section \ref{sec:intro} and defend in what follows.  One might understand the Standard View to be a view about informal proofs \emph{qua} proofs, i.e., as sources of justification for mathematical statements.\footnote{We are grateful to Harvey Lederman, Silvia De Toffoli, and an anonymous reviewer for pressing us on this point.}  It may be that informal proofs, or the activity of producing proofs, also function in other ways or for other purposes, such as by contributing to mathematical understanding, in ways that are compatible with the Standard View.  Indeed, \citet{Avigad2021} even argues that the sorts of virtues mathematicians looks for in informal proofs, the ``higher-level epistemic features of mathematical reasoning'' (7396) may actually contribute to our ability to identify when an informal proof indicates a formal counterpart, while also playing other salutory roles.  As he puts it, ``Coming to terms with the nature of mathematical justification is not at odds with understanding a wider range of mathematical values, but an integral part of the greater enterprise'' (7396). From this perspective, the position we sketched above, and the arguments we give below about epistemic value, could be seen as compatible with the Standard View.\footnote{See also \citet{Burgess+DeToffoli}, who argue that many of the disputes between advocates and critics of the Standard View are merely terminological.}  

Suppose we stipulate that there are versions of the Standard View -- perhaps even the version that Avigad or others would accept -- that are entirely compatible with our arguments.  Even so, we think it is important to emphasize that whatever formal correctness accomplishes, it is not what makes the practice of constructing and sharing proofs an integral part of mathematical activity.  In other words, it is not what makes informal proofs ``\emph{work}'' \citep[7394, emphasis ours]{Avigad2021}.  Our arguments in what follows are meant to show that formal correctness and epistemic value can come apart in ways that show that epistemic value must come from something else.  Similarly, while it may be the case that formal correctness generally coincides with mathematical correctness, in the sense that in cases where mathematicians are satisfied that they have justification for some proposition, a formal proof of that proposition could very likely be found, the formal proof is generally not what \emph{provides} the justification.  Mathematicians are not justified in believing propositions \emph{because} they have formal proofs (when they do). The dependence arrow points in the other direction.  In other words, versions of the Standard Model that are compatible with what we defend here would have very little to say about mathematical epistemology.  They would observe a correlation between formal correctness and mathematical correctness, but provide no further insight into what mathematicians do and how it works. These observations are especially important in the context of the other strand of literature we noted above, which seeks to foreground automated proof generation for future mathematics.  

\subsection{Against the Standard View}

Criticisms of the Standard View can also be found throughout the mathematical literature, at least implicitly.  These remarks most often appear in works where mathematicians reflect on their own methods and discipline, such as the now-classic remarks on ``Proof and Progress in Mathematics'' from Fields medalist William \citet{Thurston}, best known for his work on low dimensional topology and geometry; and in work by Reuben \citet{Hersh1993,Hersh1997,HershBook}.  But the philosophical literature opposing the standard view is most often traced back to \citet{Rav1999}.\footnote{Much like the present article, \cite{Rav1999} itself sits downstream from a much longer running discourse in mathematics and philosophy of mathematics.  In addition to the texts already cited, see \citet{PoincareSaus}, \citet{Polya}, and \citet{Lakatos}.}  Rav presents a thought experiment wherein we entertain a mathematical oracle that could rule on the truth or falsity of any mathematical proposition.  He suggests this would realize a Hilbertian formalist's vision of mathematical utopia---and that it would simultaneously spell doom for mathematics.  

He goes on to argue that, ``the essence of mathematics resides in inventing methods, tools, strategies and concepts for solving problems which happen to be on the current internal research agenda or suggested by some external  application. But conceptual and methodological innovations are inextricably bound to the search for and the discovery of proofs, thereby establishing links between theories, systematising knowledge, and spurring further developments'' (6).  In other words, proof is an essential part of mathematics because the activity of trying to produce, refine, and improve proofs is often the context for mathematical innovation and deeper understanding.  But the formal correctness of proof is only incidental to this role in mathematical practice, for two reasons.  First, it is very rarely within a fully formal setting that mathematicians develop significant new ideas---and of course, automated theorem provers are not the kinds of things that can have new ideas at all.  Secondly, he argues, formal proofs are some distance removed from proofs as they fruitfully contribute to mathematics, because formalizing a proof involves breaking the link to the intended semantics of the informal mathematical arguments.  Once we move to a formal setting, we are no longer talking about groups, or manifolds, or locales; instead, we are manipulating symbols according to certain rules.  For these reasons, he concludes, the Standard View misidentifies both why proofs are important and what features of proofs most contribute to that importance.

In the 25 years since Rav's essay appeared, others have also mounted attacks on the Standard View.  \citet{Detlefsen}, for instance, emphasizes that ``Mathematical proofs are not commonly formalized, either at the time they're presented or afterwards.  Neither are they generally presented in a way that makes their formalizations either apparent or routine'' (p. 17), suggesting that formalization is not a central concern for the mathematicians generating proofs, and that it is not essential to, or even especially closely related to, rigor.  Several authors focus on other properties of proofs that seem at least as important.  For instance, \citet{Easwaran} highlights ``transferability'' as a key criterion, i.e. ``that a proof must be such that a relevant 
expert will become convinced of the truth of the conclusion of the proof just
by consideration of each of the steps in the proof'' (343).\footnote{Transferability is closely related to, but a bit different from, ``shareability'' as \citet{deToffoliGFM} discusses it.  There the issue is a more general one about the intelligibility of proofs.}

\citet{Tanswell} presents a different argument, which is that in general an informal proof can be formalized in many different ways, and often in different systems, suggesting that it is a mistake to hold -- as some, such as \citet{Azzouni}, apparently do -- that an informal proof is just a way of indicating that a (particular) formal proof exists.\footnote{\citet[p. 7387]{Avigad2021} does not see a problem for the Standard View here---to the contrary, he thinks the fact that informal proofs that admit one formalization generally admit many helps explain how informal proofs can reliably convince us that any formal proof exists.}  Instead, the one-to-many relationship of informal proofs to formalizations suggests that formal proofs represent different ways of elaborating or commenting on the structure of an informal proof, and not vice versa.  Or consider \citet{Larvor}, who, responding specifically to \citet{Avigad2021}, argues that formalists often conflate two different visions of what an informal proof is: either an informal proof is itself a syntactic object (and thus just a poor version of a formal proof), or it has some semantic content.  If it is the former, then it looks like the Standard View is trivial; if it is the latter, there is a significant gap between what informal proofs actually accomplish and what a formal proof can hope to recover.  

\section{Formal Correctness is not Sufficient}\label{sec:sufficient}

We now turn to the principal arguments of this paper. We will offer two arguments that correctness is not sufficient for a mathematical proof to have epistemic value.\footnote{These arguments are similar to those of \citet{Rav1999} and others, though below we will draw connections to other mathematicians' commentary on how they see their contributions, which we think enrich the arguments.}   The first argument concerns additional necessary conditions on the assertions proved; the second concerns additional necessary conditions on the proofs themselves.

We begin with the most obvious: the fact that a proof is correct does nothing to guarantee that the proposition established by the proof is \emph{interesting}.  By interesting we mean what working mathematicians mean when they describe a result as interesting, namely, that they learn something from it that is relevant to their goals and values.\footnote{How directly relevant?  We do not mean that mathematicians only judge results to be interesting if they can immediately put them to use for their next paper.  Rather, it is that there is a shared activity, or family of activities, that mathematicians are engaged in, and they can recognize contributions to that activity.} We claim that a theorem must be interesting to mathematicians working in some area or other for it to have epistemic value.  How could a theorem fail to be interesting?  One way is for the assertion to be trivial, in the sense that no one learns anything new from the result.  Another way in which a proof could fail to be interesting is that the assertion proved concerns structures or properties that are not well-integrated with the rest of mathematics, so that even true theorems about them do not have any bearing on mathematical practice.  The proposition proved is not relevant to anything mathematicians value. And so on.  

One might be reminded, here, of the famous proof that there are no uninteresting (natural) numbers.  For suppose there were an uninteresting number.  Take the set of all such numbers.  This set contains a unique least element.  But surely \emph{that} number is interesting!  One might argue likewise about uninteresting but correct proofs: proofs are themselves rich mathematical objects, and any given proof will have features that distinguish it from others.  One might even argue that by establishing a proposition, a proof must have \emph{some} epistemic value.  The core of this response is to maintain that ``interest'' in our sense is not necessary for epistemic value.  Surely epistemic value is in the eye of the beholder, and a proof that may not be of interest to mathematicians trained to work on certain things and beholden to current fashions may nonetheless have some epistemic value in other contexts---perhaps, even, by virtue of the role it plays in computer assisted proof generation.

But this response misses the joke, and the point.  Epistemic value, here, concerns the positive epistemic contribution that proofs play in realizing the goals of actual mathematical practice.  And whatever one might take the goals of mathematics to be, surely generating sequences of sentences allowed by some proof system or other should not be counted among them.  Mathematicians do both more and less than this. The results of rote proof generation will not, in general, contribute to the epistemic activity of mathematicians.

What are mathematicians doing, beyond simply generating correct proofs?  A key goal of working mathematicians is to investigate and ultimately understand mathematical structures.\footnote{Some readers will want an account of what ``understanding'' means.  We agree that developing such an account might be interesting and would clarify aspects of our arguments, but it would also take us too far afield, and so we will simply assume that ``mathematical understanding'' is clear enough for present purposes. That being said, there is a growing literature on ``understanding'' in philosophy science \citep[e.g][and references therein]{Khalifa,deRegt}.  Some of the issues animating that literature -- for instance, does understanding require explanation? -- do not obviously bear on the sort of understanding at issue in mathematics, but others, such as whether understanding involves a kind of ability to manipulate (representations of) the thing understood, or if it is reducible to knowledge, might be more relevant.}  Very often, the most important work in this direction has nothing to do with proof, so much as with identifying definitions that both clearly and adequately capture mathematicians' ideas, and also lead to fruitful new work elaborating on what follows from those definitions in ways that provide insight into the structures mathematicians aim to describe.  

One might reasonably worry -- again -- that we are changing the subject.  The Standard View is about proofs, whereas we are now discussing definitions, proof techniques, and so on.  But this is precisely to the point.  When mathematicians speak of ``proof'', there is an ambiguity between the practice of generating proofs, which involves articulating definitions, making conjectures, developing techniques, and so on; and there are the artifacts that they produce at the end.  Our claim is that the epistemic value of proof should be located in the rich range of activities that working mathematicians engage in to produce their artifacts, not in the artifacts themselves.  Simply generating artifacts is not enough.  Indeed, we would go further: it is the practice of creating, refining, and checking the artifacts that should give advocates of the Standard View confidence that a proof can be formalized, and give other mathematicians confidence that a proof is (mathematically) correct.  The marks on the page, absent the surrounding activity leading to their generation, do not signal anything at all about formalizability (or anything else).

One aspect of this perspective is nicely captured by the Russian algebraic geometer Yuri Manin, who writes,
\begin{quote}\singlespacing Mathematicians have developed a very precise common language for saying whatever they want to say. This precision is embodied first of all in the definitions of the objects they work with, stated usually in the framework of a more or less axiomatic set (or category) theory, and in the skillful use of metalanguage (which our natural languages provide) to qualify the statements. All the other vehicles of mathematical rigor are secondary, even that of rigorous proof. In fact, barring direct mistakes, the most crucial difficulty with checking a proof lies usually in the insufficiency of definitions (or lack thereof). In plain words, we are more deeply troubled when we wonder what the author wants to say than when we do not quite see whether what he or she is saying is correct. The flaws in the argument in a strictly defined environment are quite detectable. Good mathematics might well be written down at a stage when proofs are incomplete or missing, but informed guesses can already form a fascinating system: outstanding instances are A. Weil’s conjectures and Langlands’s program, but there are many examples on a lesser scale. \citep[p. 166]{Manin}\end{quote}
Definitions are where mathematicians identify what they find interesting.\footnote{At least some leading early formalizers were also vividly aware of this, e.g. \citet[pp. 493-494]{Huntington}.}  Without the right definitions, proofs are unlikely to be of epistemic value.

Of course, what follows from good definitions is also important, and that is established by proofs. We do not mean to say that proofs do not matter.  But we would argue that what makes proofs epistemically valuable is that they show not only \emph{that} some assertion or other is true, but also \emph{why} it is true, by showing how it follows, using accepted techniques, from prior results that one already understands and accepts.\footnote{This is not quite to say that proofs are mathematical explanations---or to take any particular view of what mathematical explanation amounts to \citep[see][\S 2]{Mancosu}.  Our point is more mundane: proofs involve drawing connections between different bits of mathematics in a way that is generally illuminating.}  Now, one can argue that some proofs are more explanatory than others, or that some theorems of interest may not admit of proofs that are clear and digestible to humans, and yet merely knowing if they are true is of value because of the role they play in proving yet other claims.  But none of this contravenes the fact that all else being equal, simple, perspicuous, and comprensible proofs are of special value;\footnote{See \citet{BillQuadratic} for a case study of this particular sort of value.} and proofs that establish a theorem without being understandable to \emph{any} human mathematician do not contribute to the activity mathematicians are engaged in.  A necessary condition for a proof to have epistemic value is that it contributes to mathematicians' understanding of the structures they seek to study---specifically, mathematical structures that are believed to be well-integrated with the rest of mathematics, and ultimately, with the natural and social worlds.

Taken together, these arguments are meant to establish that formal correctness (or rigour) is not sufficient because formally correct proofs do not necessarily establish interesting propositions in ways that contribute to mathematicians' understanding of their subject matter.  Again, these arguments are not intended to be novel: to the contrary, our goal here is to isolate the core of ideas that are ubiquitous among mathematicians, and frequently expressed by leading figures.  Take, for instance, Thurston: 
\begin{quote}\singlespacing
[W]hat [mathematicians] are doing is finding ways for people to understand and think about mathematics.

The rapid advance of computers has helped dramatize this point, because computers and people are very different. For instance, when Appel and Haken completed a proof of the 4-color map 
theorem using a massive automatic computation, it evoked much controversy. I interpret the controversy as having little to do with doubt people had as to the veracity of the theorem or the correctness of the proof. Rather, it reflected a continuing desire for human understanding of a proof, in addition to knowledge that the theorem is true. 

On a more everyday level, it is common for people first starting to grapple with computers to make large-scale computations of things they might have done on a smaller scale by hand. They might print out a table of the first 10,000 primes, only to find that their printout isn't something they really wanted after all. They discover by this kind of experience that what they really want is usually not some collection of ``answers''—what they want is understanding.  \citep[p. 162]{Thurston}
\end{quote}
Thurston goes on to argue that an ``...emphasis on theorem-credits has a negative effect on mathematical progress'' \citep[p. 172]{Thurston} because it obscures both the importance of understanding as an output of mathematical practice and the importance of understanding acquired from others in producing proofs in the first place.

As with the quote from Manin, this passage from Thurston is just one well-known example of a prominent mathematician expressing this sort of view.  But we need not multiply examples to make the point, since we find that it is widely recognized even among advocates for the Standard View. 
 Take \citet[\S\S 3-5]{Avigad2022}, for instance, who develops an account of mathematical understanding and mathematical depth, and acknowledges both that mathematicians place great value on understanding and also that formal proofs, or at least, some prominent computer generated proofs, apparently do poorly on delivering understanding.  (He also acknowledges the importance of good definitions.)  But he goes on to argue that perhaps there is a different kind of understanding that at least some formal proofs can provide, even if it is not quite the kind of understanding mathematicians have typically sought out.  We are happy to grant that this can happen.  But it is not an argument that formal correctness is sufficient for a proof to have epistemic value.  At best, it establishes that there may be multiple ways in which proofs can have epistemic value.

\section{Formal Correctness is not Necessary}

We now turn to showing that formal correctness is not necessary for a mathematical proof to have epistemic value.  This direction is more subtle, and we find fewer precursors in the literature.  Consider that on a common strategy to establish that one property is not necessary for another, one wishes to exhibit, or at least show the existence of, something that instantiates the second property but not the first.  In the present case, that would mean presenting a proof that has epistemic value but which is not formalizable.  Of course, this is possible.  As even the defenders of the Standard View freely admit, real, published mathematical proofs often contain errors, and as such they are not formalizable as-is.  And yet they may nonetheless have considerable epistemic value: the errors may be trivial, so that when corrected, the proof is formalizable; or the proof may present a successful proof strategy even though it does not accomplish its goals as implemented; or the erroneous proof may have a gap, identifying which advances the field in some way; or the proof may introduce new techniques that turn out to be fruitful elsewhere.  

One might justly complain about this argument, however.  First, proofs with simple errors fail to be formalizable in a trivial way (and, indeed, \citet{Avigad2021} argues that standard mathematical methods are effective in part because they are robust against small errors, so that erroneous proofs are effectively formalizable).  Meanwhile, to cite ``proofs'' that have more serious errors, but which contribute to mathematical knowledge in other ways, is arguably  question-begging.  ``Proof'', here, should be understood as a success term.  An argument with significant errors, ones that are not easily corrected, is \emph{not a proof} in the relevant sense, even if it is typeset as a proof and published in a mathematical journal.  At best one can show that there are (fallacious) arguments, or proof-attempts, or simil-proofs that have epistemic value.  (We will do exactly this below.)  And of course, it is unsurprising that fallacious arguments cannot be formalized as formal proofs in a consistent system.  So what one would presumably need to do to make the argument go through would be to find an example of a proof that is widely accepted by the mathematical community as error-free, convincing, and otherwise laudable, but which is nonetheless not formalizable.

We will not attempt to give an example of this sort.  We do not think it is possible.  Even if one offered a compelling candidate, we are skeptical that a convincing argument could be given that such a proof is not ``formalizable''.  The reason is the essentially modal character of the demand: the question is whether there is \emph{some} formal system or other such that the proof could be formalized as a (valid) proof in that system.\footnote{One might think it would suffice to give an accepted proof that cannot be formalized in \emph{standard} systems.  We return to how to think about that possibility at the end of section \ref{sec:anti-logicism}.}  This requirement is largely unconstrained.  Moreover, given the nature of mathematical practice, and the standards for systematicity and rigor that mathematicians demand, one would expect that the arguments of any widely accepted proof could be made sufficiently explicit in some formal system or other---and if it turned out that existing formal systems were not up to the job, it would be an interesting project in mathematical logic to try to identify a new proof system in which the proof \emph{could} be formalized.

One might think that conceding this point is to concede that correctness is necessary after all.  But we wish to point to three senses in which it is not.  Two of these are subtle, and we address them in the next section.  But the first is simple: if we drop the strongly modal version of the requirement, then formalizability is clearly not necessary, by the very same reasoning. In other words: there is no formal system such that a given proof must be formalizable \emph{in that system} for it to be epistemically valuable or otherwise mathematically acceptable.  This claim is hardly new.  At least on one reading \citep{vonNeumannHilbert}, Hilbert's program in foundations of mathematics was to show that any correct proof could be formalized in first-order arithmetic (or, more generously, some fixed set theory)---which, in turn, he hoped, could itself be shown to be consistent \emph{within that same system}.  But of course, this goal foundered on the shoals of G\"odel's incompleteness theorems, because any system in which Peano arithmetic (say) could be shown to be consistent would have to be strictly stronger than Peano arithmetic.  Thus, when \citet{Gentzen} later showed that Peano arithmetic is consistent, his proof assumes transfinite induction up to $\epsilon_0$.  And if one wished to show that Gentzen's system is consistent, one would need a stronger system still.  And yet, the most committed advocates of formal correctness as a guiding norm in mathematics would surely want to concede that these sorts of arguments in the foundations of mathematics are epistemically valuable.\footnote{To be sure: we are not claiming that aiming to prove consistency results is deeply connected to the Standard View.  Consistency results in mathematical logic should be seen as mathematical theorems like any other.  Our point is that you cannot fix a single formal system and insist that epistemic value is identified with formalizability in that system.}

One can see a different dynamic at play in the classic debate over intuitionism in the foundations of mathematics, though it leads to the same basic moral.  For early intuitionists like Brouwer or Weyl, a proof could be epistemically valuable only if it was \emph{constructive}, i.e., if existence proofs proceed by showing how to construct an instance of the things whose existence is claimed.  This idea has led to several proposed systems of intuitionist logic, including the Brouwer-Heyting-Kolmogorov interpretation of constructive logic and Martin-L\"of's intuitionistic type theory.  Various authors have argued that mathematics must be formalizable in some such system in order to be epistemically valuable in the fullest sense, on the grounds that non-constructive proofs are unreliable.  But such approaches have never gained widespread acceptance within mathematics, because intuitionistic approaches tend to be weaker than the methods of standard, informal mathematics.  They limit what mathematicians can do.  And when pressed, working mathematicians prefer to follow fruitful mathematical developments where they lead, rather than to constrain themselves to particular formal systems.

Towards the end of his life, reflecting on the debates in the foundations of mathematics on which he had cut his teeth from the perspective of someone who had shifted towards more applied questions, von Neumann summarized his attitude towards the earlier debates that we have been discussing as follows.
\begin{quote}\singlespacing
    In my own experience, on two other occasions in the early twentieth century, there were very serious substantive discussions as to what the fundamental principles of mathematics are; as to whether a large chapter of mathematics is really logically binding or not.  And in the nineteen-tens and -twenties a critique of these questions made it apparent, that it was not at all clear exactly what one means by absolute rigor, and specifically, whether one should limit oneself to use only those parts of mathematics which nobody questioned.  Thus, remarkably enough, in a large fraction of mathematics there actually existed differences of opinion!  Some mathematicians said that one need not question any part of what is in fact being used.  There was also a body of opinion, that one should not use more than what the most exacting critics had approved.  However, there was a further, large body of mathematicians, who felt that while there was some point in questioning certain areas of mathematics, it was all right to use them.  This group was quite ready to accept something like this: Those portions of mathematics which had been questioned and which had been clearly useful, specifically for the internal use of the fraternity---in other words, when very beautiful theories could be obtained in those areas---that those were after all at least as sound as, and probably somewhat sounder than, the constructions of theoretical physics.  And after all, theoretical physics was all right; so why shouldn't such an area, which had possibly even served theoretical physics even though it did not live up to 100 per cent of the mathematical idea of rigor, why shouldn't this be a legitimate area in mathematics; and why shouldn't it be pursued? This may sound off, as well as a bad debasement of standards, but it was believed in by a large group of people for whom I have some sympathy, for I'm one of them. \citep[pp. 480-1]{vonNeumannSociety}
\end{quote}
When confronted with a choice between adhering to some fixed standard of rigor -- of correctness -- and valorous mathematical theories -- beautiful ones, or deep ones, or fruitful ones -- a significant portion of mathematicians will always prefer the latter. We contend they are right to do so---and for precisely this reason, formalizability in any particular system cannot possibly be deemed necessary for epistemic value in mathematics.  

\section{Fruitful Errors}

In the previous section, we made the argument that formalizability in any particular system cannot be necessary for mathematics.  We concluded with a passage from von Neumann pointing out that in cases where fruitful mathematics was found to run afoul of some formal system or other, mathematicians tended to reject the formal system and keep the fruitful mathematics.  But perhaps this argument targets a too-strong version of the Standard View, and formalizability should be understood more broadly.  In this section, we wish to reflect further on a suggestion we find in the passage from von Neumann -- taking, perhaps, some liberties of elaboration --  to the effect that what ultimately grounds the epistemic value of some area of mathematics, or even some particular proof, is how that piece of mathematics fits into a broader network of mathematical and scientific activities.  

We have used the term ``fruitful'' several times here.  What we have in mind is that fruitful mathematics is mathematics that leads to yet more of the same sort of activity, to new ideas and applications in physics or other mathematical sciences, or even just to yet more mathematics.\footnote{Fruitfulness, here, is closely connected to what \citet{MaddyDtA} calls ``depth'' \citep[see also][]{Maddy+etal}.}  Mathematics that is not rigorous but which leads to successful applications can very likely be made rigorous, though doing so can itself be a challenging (and rich) task;\footnote{Calculus provides a prime historical example, having only been made rigorous to modern sensibilities two centuries after its introduction; another example is the Dirac delta function.  The path-integral formalism of quantum field theory will very likely follow a similar trajectory.} but mathematics that is formally correct but not guided by fruitful application is not at all certain to ever find fruitful application---for reasons we have already argued for above.

This observation leads us to the second sense in which correctness is not necessary for mathematics, alluded to above: as a purely descriptive matter, within ordinary mathematical practice, correctness does not need to be established in order to secure the broad acceptance of a new result.  Clearly, as \citet{Detlefsen} emphasizes in this connection, it is not necessary to first formalize a new mathematical result before publishing it in even the most prestigious mathematics journals.  This, of course, is the observation that has driven much of the debate over the Standard View; and which has led \cite{Avigad2021} to argue that the norms and strategies of informal proof can be seen as guides towards formalization without actually securing formalizability.  But in fact, something deeper is true.  Cases of mathematical \emph{disagreement} are almost never resolved by formalizing a difficult and controversial argument.  When a major conjecture is settled using new techniques, for instance, and the mathematical community does not immediately accept it, one does not find mathematicians on each side of the dispute rushing to produce formalizations or counter formalizations of the disputed proof.  New techniques are not certified in this way.  

What mathematicians \emph{do} care about is ensuring that the new piece of mathematics is sound, in the sense that it preserves truth about whatever the target structure under consideration is.  They want to be sure that counterexamples to the new result cannot be constructed using standard methods, that they will not inadvertently introduce contradictions by assuming the result holds, and that future applications of the theory will be reliable.  In principle, formalization could perhaps help with this.  But it is not the most straightforward, perspicuous, or common way of doing it.  This is because translating informal proofs into formal ones is not, in general, easy; it does not generally get to the heart of the subject matter; and \emph{failures} to formalize in some particular system are not necessarily probative, for reasons already discussed.  And so correctness is not necessary in the colloquial sense,  since it is not required or needed to convince mathematicians of the reliability of a disputed result.

Finally, we turn to a third argument.  We have already conceded that ``proof'' is a success term, and that plausibly a proof that cannot be formalized, in any system at all, should not be deemed a proof.  We have also argued that there are still senses in which formalizability is not necessary.  But what if we forget about proofs, per se, and consider instead the more general category of strong proof candidates -- roughly, simil-proofs -- including serious attempts at proof that may even appear in the published literature, but which turn out to be incorrect.  There are many examples, we claim, of proof candidates that ultimately turned out to be incorrect---but which had epistemic value, as witnessed by their impact on subsequent mathematics.  We will describe two such examples here.  Our claim is that these non-proofs had as much, or more, epistemic value as many important proofs.  In fact, we will argue, although they were erroneous, they also helped mathematicians more clearly grasp the shape of novel and previously unsuspected phenomena of fundamental importance to later mathematics.


In 1895, Poincar\'e released the first installment of his monumental paper {\em Analysis Situs}.  His stated goal in this paper was to launch a new field of mathematics, now known as ``algebraic topology''.  To say that he succeeded is to do scant justice: algebraic topology grew throughout the twentieth century and into the present to be one of the major fields of mathematical research, with  influences and cross-pollination across mathematics, and across the sciences more generally.  Today it finds  significant applications in condensed matter physics, quantum chemistry, and contemporary data science, just to name a few.  {\em Analysis Situs} is also riddled with errors, which led Poincar\'e to produce five additional ``supplements''.  Even these did not suffice to fully settle the issues, neither to today's sensibilities nor to those of his contemporaries. 

It is hard to square the success and impact of these papers with the formalist view of mathematics.  Standard formalist responses have alternated between treating them as anathema and ignoring them altogether.\footnote{For instance, ``For Bourbaki, Poincar\'{e} was the devil incarnate'' \citep[p. 280]{Mandelbrot}.  The quote is also attributed to Stone \citep[p. 145]{MacHale}.} Our goal here is not to rehash this discussion, which others have engaged with before \citep[e.g.][]{PoinStill,Sarkaria}.  Rather, we want to focus on the fecundity of Poincar\'e's mistakes, most notably his false assertion that every homology 3-sphere is homeomorphic to a 3-sphere, which led to his discovery of his dodecahedral space and the formulation of the ``Poincar\'e conjecture''.\footnote{For a careful discussion, see the translator's introduction to Stillwell's translation of {\em Analysis Situs} \citep{PoinStill}.}  This has been one of the most fertile mathematical conjectures of the past century and a half, leading to Smale's work on $h$-cobordism (and the proof of the conjecture in dimension $\ge 5$) \citep{Smale}, Freedman's work on the disc-embedding theorem and the classification of simply connected (topological) 4-manifolds \citep{Freedman}, and Thurston's geometrization conjecture \citep{ThurstonGeom}/Perelman's theorem classifying the (piece-wise) uniform geometries on closed 3-manifolds \citep{Perelman1,Perelman3,Perelman2}. For those keeping score at home,  these results alone comprise three Fields Medals awarded, one Fields medal declined and one (\$1,000,000) Millenium Prize declined (the only one which could have been awarded to date). Notably, Poincar\'e's errors are only the first in a series of major errors which have given rise to large swaths of contemporary algebraic topology and algebraic geometry, including:
\begin{enumerate}
    \item Pontrjagin's famous incorrect computation of the stable homotopy groups of spheres in \cite{Pontrjagin38}, leading to the ``Pontrjagin-Thom'' construction \citep{Pontrjagin47,Pontrjagin55,Thom}, Thom's work on cobordism theory \citep{Thom}, and the formulation in \citep{Kervaire} of the ``Kervaire invariant one'' problem (solved by \citet{HillHopkinsRavenel}, except for one open case whose solution has recently been announced in \cite{LinWangZhu}).
    \item Lefschetz's work on the topology of algebraic varieties, which led Hodge to initiate what we now call ``Hodge theory''\footnote{Hodge theory itself was also initially beset by errors. As Atiyah wrote in in Hodge's 1976 biographical entry for the Royal Society, ``In retrospect it is clear that the technical difficulties in the existence theorem did not really require any significant new ideas, but merely a careful extension of classical methods. The real novelty, which was Hodge’s major contribution, was in the conception of harmonic integrals and their relevance to algebraic geometry. This triumph of concept over technique is reminiscent of a similar episode in the work of Hodge’s great predecessor Bernhard Riemann.''} and the ``Hodge Conjecture'' (an open Millenium Prize problem). A classic quip, repeated by Griffiths in his biography of Lefschetz 
    \citep[p. 289]{NAP2037}, asserts that ``Lefschetz never stated a false theorem nor gave a correct proof.''  Notably, one of the core theorems in question is a refinement of Poincar\'e duality, the site of a different major error in {\em Analysis Situs} and the one discussed by \citet{McLarty}.
\end{enumerate}
And this is to name only two of the most consequential and far-reaching errors! 

For another example, consider Frege's monumental effort to formalize mathematics \citep{FregeGrund1, FregeGrund2}, which is frequently credited as the origin of the current formalization frameworks. On the eve of publication of the final volume of his attempt to derive all of mathematics from logic applied to a few self-evident axioms, Frege received a letter from Russell showing that the naive set theory underpinning Frege's project contained unavoidable contradictions. In an appendix written as the book was going to press, Frege described the experience as follows: 
\begin{quote}\singlespacing
    Hardly anything more unfortunate can befall a scientific writer than to have one of the foundations of his edifice shaken after the work is finished. This was the position I was placed in by a letter of Mr. Bertrand Russell, just when the printing of this volume was nearing its completion.\footnote{Appendix of \cite{FregeGrund2}, in \cite[p. 279]{Beaney}, translation by Michael Beaney.}
\end{quote}
From one view, Russell's paradox mortally wounded Frege's project, with the final nail in the coffin delivered several decades later by G\"odel's incompleteness theorems.  However, as a historical matter, one could hardly hope for a more influential or productive mistake. The drama of Frege's ``failure'' provided a significant impetus and status-boost to research in axiomatic set theory and mathematical foundations, leading to several decades of work and attracting leading mathematicians including Zermelo, Hilbert, Ackerman, von Neumann, G\"odel, Bernays, and Tarski, just to name a few. Even more, the error itself was and is a major discovery: every sufficiently expressive\footnote{``sufficiently expressive''=``sufficient to express integer arithmetic''.} finite formal system admits self-referentiality as an inescapable feature, and this radically limits the behavior of such systems. As a mathematical idea, this continues to bear fruit, from G\"odel's second incompleteness theorems showing that no {\em finite} formal system is capable of satisfying Hilbert's criteria for an axiomatic framework for mathematics, to Turing's method for showing undecidability in theoretical computation, and then the many problems of mathematical and scientific interest for which no general algorithmic solution can exist.\footnote{At least in the Hilbert-Church-Turing paradigm of algorithm.} 

We need not multiply examples further---though we insist that cases like these are not wildly unusual, except in their overall importance to the development of mathematics.  These already suffice to establish the point that contributions to mathematics, including would-be proofs, need not be formalizable or even correct to have epistemic value.  Now, one might object that whatever value these examples of mathematical output might have, it is not the value that proofs have---or at least, not the value that proofs have \emph{as proofs}, which is specifically the value of \emph{justifying} belief in a mathematical proposition, because they are not truly proofs.  But, for reasons we have already given, we are skeptical that what makes proofs valuable to mathematicians is really that they justify propositions, or at least, they do not do so by virtue of being formalizable.  

Instead, we have argued, the value of proof comes primarily from how they permit mathematicians to more clearly reason about and understand mathematical concepts.  Of course they do contribute to mathematicians' beliefs about what propositions are true and false, but we suggest they do so through the role they play in a more general practice.  Mathematicians come to justify mathematical propositions through a range of methods: trying to find counterexamples, applying known techniques and developing new proof techniques, drawing connections with other areas of mathematics, showing how those propositions depend on other propositions, as lemmas and conjectures, and so on.  This process often results in an (informal) proof, but it is the whole process that ultimately provides justification.  And this is what makes these examples so important, because they seem to show how \emph{incorrect} arguments can also contribute, in just the same way, to that more general practice of exploring mathematical terrain and coming to better understand mathematical concepts. If one understands proof as functioning in this way, then these errors are fruitful for coming to learn what is true in just the same way, and sometimes to a greater degree, than correct arguments.  

One might worry at this point that we have shown too much, and that our arguments in fact show that even mathematical correctness -- or rigor, as ordinarily understood in mathematics -- is not necessary for epistemic value.  After all, the problem with Poincar\'e's claims was not just that they were not formalizable; it was that they were simply \emph{wrong}, by his own lights and those of the mathematical community.  And yet, we claim, they still had epistemic value---and even that they contributed to justification.  This is a serious worry.  But it is important to clarify that we do not think incorrect proofs can provide justification for anything in and of themselves.  Ultimately, we need (mathematically) correct arguments to have satisfactory justification.  Our point is rather that incorrect proofs contribute to the process of developing and assessing those arguments, and that epistemic value, including justification, should be located with the process rather than the final result.\footnote{Recall the discussion in section \ref{sec:sufficient}.}  

Perhaps more importantly, mathematicians' own assessment of the rigor and mathematical correctness of their work and that of their colleagues is dynamic, and very often backwards looking.  And so, while it is ultimately necessary, it is not something that can be immediately assessed.  In key examples from the history of mathematics where important errors went overlooked -- as (at least initially) with Poincar\'e, but also other examples, such as Cauchy's work on uniform continuity -- the relevant mathematicians believed their work \emph{was} rigorous, and it was often assessed as rigorous by the community.  Indeed, many of the formal methods we have today were developed only in response to later mathematicians finding problems in earlier mathematicians' work, and trying to identify and avoid those problems.  Note however that this clarification did not take place through some rote process, but rather through mathematicians more deeply understanding the core phenomena and developing new methods, including but not solely formal ones, that could help avoid those errors.  And so, even here, we would say that while mathematical correctness is ultimately necessary for proofs to provide justification, assessing mathematical correctness is itself part of the process we intend to emphasize---a process that often also includes making, and identifying, errors.

Now, as we have acknowledged above, one might also object that what we have said is compatible with the Standard View, insofar as advocates for the Standard View may feel that proofs play \emph{many} roles, and formalizability is essential only for one of those roles---one that such advocates could even concede is relatively minor compared to other roles played by proof. Fair enough.  But even so, we think the Standard View is misguided, because it reverses the relationship between formalizability and mathematical correctness or rigor, for reasons we discuss in the next section.

\section{The Relationship Between Logic and Mathematics}\label{sec:anti-logicism}

Our discussion to this point has concerned ``formal correctness'', that is, correctness as determined by (valid) formalizability in some formal system.  We have argued that formal correctness is neither necessary nor sufficient for epistemic value in mathematics.  This raises the question of whether formal correctness bears any relationship to the epistemic value of mathematical proof at all.  We suggest that the answer is ``yes'', but the situation is subtle.  Our answer goes via a different notion of correctness, which might be called ``mathematical correctness'' (or ``informal correctness'').  ``Mathematical correctness'', here, is the standard of correctness for mathematical proof accepted by the mathematical community, on reflection. Somewhat more precisely, we take a lead from Peircean pragmatism and define mathematical correctness as ``correctness at the end of inquiry'', i.e., after indefinite scrutiny and in light of all subsequent mathematical developments.\footnote{For those worried about the apparent requirement of infinite time here, we suggest hearing this ``weakly''; we are not aware of an instance where a span of more than a century of active research was required to detect an error in a central, widely credited result.  That said, one might worry about whether \emph{interest} would shift over time, so that in-fact incorrect results are never identified as such because the community moves on to other issues.  For present purposes, we do not worry about this sort of case.}  Thus, results in early 19th century analysis may be viewed as mathematically incorrect even though they were accepted by the community at the time.  On the other hand, we do not want to suggest that the sociological whims of future mathematicians can retroactively dictate what makes good mathematics.  One way to think about this worry is to say (somewhat speculatively) that Cauchy, Dirichlet, and their cohort would have accepted later arguments that aspects of their work were not mathematically correct, and in that sense the standards do not change.  

Advocates for the Standard View would maintain that formal correctness is equivalent to mathematical correctness, or at least, necessary and sufficient for mathematical correctness.  Indeed, there is a long tradition in philosophy of logic and mathematics, originating, arguably, with Leibniz's \emph{characteristica universalis}, of conceiving of logic as providing a normative metaphysical and perhaps epistemic foundation for mathematics---or, in its most totalizing form, all human knowledge. This picture is one on which the rules for manipulating formal syntax are taken to be fundamental laws of reason, and various logical axioms are necessary truths that are prior to other forms of knowledge.  From this perspective, logic is fundamental.  Mathematics, meanwhile, is to be reduced to, or at least its truth is to be grounded in, logic.  Arithmetic is true insofar as it can be reduced to a logical system; perhaps even our knowledge of arithmetic ultimately comes from its relationship to the secure truths of logic.  

We endorse a different view of the relationship between logic and mathematics---one we think has a folk status in the literature, but which is also directly endorsed by some authors, such as \citet{MaddyEIT}.\footnote{The view is also anticipated by \citet{BurgessProofs}, though somewhat less explicitly, and we suspect versions of it can be seen in other work as well.}  On our view, mathematical logic is best seen as a mathematical theory, or family of mathematical theories, of mathematical practice.  It bears a relationship to ordinary mathematics analogous to the relationship that a mathematical theory in science bears to some real-world phenomenon.  It captures some aspects of mathematics in precise mathematical terms.  It allows us to prove theorems about what mathematics (so modeled) can and cannot accomplish, and it allows us to reason about the relationship between sentences, truth conditions, and structures that may instantiate what sentences assert.  Among other things, it provides a mathematical framework for thinking about how mathematicians reason, and it offers a precise account of mathematical correctness---namely, formal correctness.  

From this point of view, mathematical logic is a valuable tool for thinking about mathematics.  But to take formal correctness as the standard against which mathematical correctness is evaluated, or worse as equivalent to mathematical correctness or as the feature that makes proofs epistemically valuable as proofs, is to put the cart before the horse.  It is analogous to saying that the reason the earth follows its orbit around the sun is that our best theory of gravitation says it ought to.\footnote{Of course, it may well be that the central equations or laws of gravitational theory assert relationships that truly obtain in the world, and that those relationships are the reason that the earth follows its orbit around the sun.  The point is that our theories attempt to describe the regularities that obtain in the world, rather than that the regularities in the world are somehow obligated to conform with our theories.}  Or perhaps better, to say that a proof is mathematically incorrect because it is not formalizable would be like saying that Mercury has erred in some way because its perihelion advances a few dozen arcseconds each century relative to its Newtonian orbit.

We do not claim that the originators of mathematical logic thought of their work in the way we have just described.  (Perhaps some did.\footnote{For example, \citet[p. 74]{Gentzen} indicated that his goal was to develop ``a formalism that reflects as accurately as possible the actual logical reasoning involved in mathematical proofs''.  Some philosophers read this ``as meaning that all valid mathematical inferences should be instances of the logical rules of natural deduction, or closely related to them'' \citep[p. 480]{Tatton-Brown}, but we find that strained.  The plain meaning of the text is that mathematical inference comes first, and the formalism of natural deduction was developed in its image.}  We are not sure.  Others certainly did not.)  Instead, we think of this view as a rational reconstruction, one that makes the best sense of the accomplishments of mathematical logic while also keeping track of its correct normative role in mathematical practice.  

So why hold the view?  The first reason is that mathematical logic clearly \emph{is} a mathematical theory of mathematical practice.  Classical first-order logic captures many aspects of mathematical reasoning.  First-order model theory does describe a relationship between logical syntax and argumentation.  Second-order logic captures common reasoning in fields like point-set topology that typically use higher-order quantification.  Intuitionistic logic captures the reasoning of mathematicians who reject certain proof methods allowed in first-order classical logic.  And so on.  All of these are mathematical theories.  So they are mathematical theories of (aspects of) mathematical practice.  We infer from this that formal correctness, relative to a given system, is a precise way of capturing mathematical correctness within that system.

So far so good.  But the important question is whether mathematical logic, and by extension formal correctness, is also \emph{more} than this.  In other words, one might argue that although mathematical logic is a mathematical theory of mathematical practice, it has a saliently different relationship to its subject matter than other mathematical theories of some natural phenomenon.  We accept it does have some important differences, such as that it is self-referential.  Insofar as mathematical logic is itself a part of mathematical practice, we can apply it to itself.  But this is not a salient difference, because self-referentiality does not imply that mathematical logic has a normative status with regard to its subject matter, whereby failures of descriptive accuracy imply that the target is erroneous.

One possible basis for a normative status of this kind would be if mathematical logic were well-conceived as an ideal form of mathematical reasoning, something like a perfection, precisification, or elaboration of ordinary mathematical practice.  In that case, one might argue that mathematical logic is mathematical reasoning with certain failure modes removed, or that if a proof cannot be formalized, it fails some higher, more pure or ideal standard.  But we reject this.  Yes, mathematical logic describes a certain kind of ideal of reasoning -- or multiple ideals, associated with different logics -- but that is because it is a highly \emph{idealized} theory of mathematical reasoning.  By this, we mean that various bits of mathematical logic describe more or less simplified versions of mathematical reasoning, with domains of application in which they are more or less descriptively accurate, and other domains where they break down.  

In fact, we have already seen this.  First-order logic has a number of nice properties that capture deductive reasoning, it has a clear semantics, and one can prove fundamental theorems -- of soundness and completeness -- that link syntax and semantics.  But standard methods in some areas of mathematics involve, and arguably require, second-order quantification.\footnote{For a full-throated defense of the importance of second-order logic for the foundations of mathematics, see \citet{Shapiro}.  One might worry that second-order logic, or second-order quantification, is not really what is needed in most cases where first-order logic fails, and that working mathematicians rarely explicitly quantify over sentences.  Fair enough---but this is a dispute about what is the best alternative in cases where first-order logic breaks down.  Another candidate, at least in some cases, is to build the relevant structures in set theory, though doing so hardly captures all aspects of the relevant mathematical practice.}  This signals a breakdown of first-order logic as an adequate model for mathematical practice.  Of course, it does not mean that these methods cannot be formalized---one just needs to move to second-order logic.  But second-order logic, while apparently expressive enough to capture classical mathematics, does not have the same nice semantic or proof theoretic properties as first-order logic.  It is famously a consequence of G\"odel's incompleteness theorems that under standard semantics, there is no (finitary) proof system that satisfies all three of soundness, completeness, and effective proof checking.  One can adopt alternative semantics (most notably, Henkin semantics), relative to which the standard deductive systems for first-order logic become sound, complete, and effectively checkable, but Henkin semantics reinterprets second-order quantification in a way that effectively reduces its expressive power, limiting its ability to faithfully model mathematical practice in the wild.  So while second-order logic improves on first-order logic in some regards, it breaks down under other circumstances.  Neither seems to perfectly capture mathematical practice, though both are helpful in certain domains.

Something similar can be said for mathematical arguments that move between proofs within a theory to proofs at the semantic or metamathematical level.  As \citet{Rav1999} points out, there are many examples in mathematics -- he highlights group theory in particular -- where a first-order theory is available, and for some purposes, textbook theorems can be seen as informal arguments that would admit formalization in a sound first-order deductive system.  But often, powerful theorems of group theory require one to manipulate groups using the full resources of set theory, that is, to step outside of the first-order theory of groups and prove things externally about structures that satisfy those axioms.  Or consider that for standard results in synthetic differential geometry, one needs to move to an intuitionistic logic and drop the law of the excluded middle.  In that case, there are good reasons to think the results are mathematically correct, because the geometric theory is dual to an algebraic theory that does not require any changes in logic.  But a flat-footed attempt to recover those results as formally correct leads to contradictions.

Of course, none of these arguments imply that formal correctness, relative to some proof system, is not a valuable tool for assessing mathematical correctness.  In many cases, showing that a proof is not formally correct, or could not be formally correct, relative to some system is probative in assessing whether it is mathematically correct, since the failure of formal correctness can reveal unnoticed errors.  But we should think of this as a case of reasoning about some subject matter by using a simplified model that we expect to represent salient features for the case at hand.  It is analogous to an engineer who uses mathematical models of material strain in the course of designing a bridge.  The models are not the \textit{reason} the bridge will stand or fall, and they can even be wrong under some circumstances.  But they are likely a useful guide nonetheless, and will point to the places where the bridge design may be flawed.

We conclude this section with a final argument, which is somewhat speculative but which seems to us to accurately reflect both the sociology of mathematics and the relationship between logic and mathematics.  Perhaps it is a Rorschach test of sorts.  Suppose that there were a result that could not be formalized in \emph{any} known logical system, but which nonetheless the mathematical community came to accept as mathematically correct.  This may seem implausible---though we would argue the reason it seems implausible is that mathematical logic is a mature discipline, and there are many systems available already.  But suppose it happens.  What should we make of that situation?  Advocates of the Standard View would presumably say that the mathematical community has gone off the rails, they have accepted an incorrect result, where incorrectness is determined by failures of formalizability.  (After all, it is hard to see in what sense the informal proof has indicated, much less encoded, a formal proof.)   But that is not what we would expect mathematicians or logicians to actually do.  Instead, we think such a situation would spark a lot of excitement among logicians, and important new work trying to find a \emph{new} formal system in which the result could be formalized.  One would try to understand how the result works and why it might be taken to be correct.  In other words, if existing mathematical theories of mathematical practice failed in some domain, the correct and likely response would be to develop new mathematical theories of that domain.

\section{The Automation of Mathematics Revisited}

As we said in the introduction, one motivation for our arguments is to reflect on what automated theorem provers and AI mathematicians can contribute to mathematics---and to push back against the idea that automating mathematics is certain, or even likely, to significantly advance mathematics.  We now return to this motivating theme.

As we noted in the introduction, if one thinks that formal correctness is the guiding virtue of mathematical practice, if one believes it is sufficient for epistemic value, then it would follow that a computational system that could reliably generate correct proofs would generate epistemic value.  This would be all the more true if the computational system could reliably generate proofs faster than humans could generate them; or if it could generate proofs of theorems that humans have not been able to, or even in principle could not, prove.  If, on the other hand, correctness is neither necessary nor sufficient for a proof to have epistemic value, it is much less clear what reliably generating correct proofs can accomplish.  It would seem that mathematicians need more than this.  Where does that leave us?

The first thing to say is that there is a sharp distinction to draw between computational proof generation and proof assistance.  Using computers to help prove theorems that have been stated by mathematicians, which assert things that mathematicians have antecedent interest in establishing, is surely valuable.  Proof assistants can check potentially problematic lines of reasoning and provide some insight into their status; in principle, they can also fill gaps and play other kinds of supporting roles.  None of this is problematic, as long as the crucial role that human mathematicians play in posing questions and extracting new understanding from answering them remains in focus.  For that to be the case, the inputs -- the conjectures stated, the definitions used, etc. -- must come from mathematicians who are seeking to extend the body of existing mathematical knowledge, with full integration with the rest of that body of knowledge; and the outputs, namely, proofs of theorems, must be such that they can be read and understood by human mathematicians. Even more, these outputs must provide the kind of insight into the structures under consideration, and the reasons theorems are true, that other, non-assisted work does.  In other words, proof assistants can surely play a collaborative role in producing new mathematical research, as long as the overall character of the resulting inquiry is not adversely affected.

More troubling, though, is the prospect that generative AI, such as high-end large language models, could \emph{simulate} the sorts of inputs and output just described.  That is, a system trained on mathematical definitions could surely state definitions that are similar in form to ones mathematicians would state; and they could pose conjectures with the right sort of form; and then automated theorem provers could generate proofs of them.  One could even imagine (and there are ongoing efforts to build) AI systems that could take the output of a formal proof assistant / generator and rewrite it in the style of a working mathematician.  In some ways this would be similar to generating prose in the style of various authors, or generating working computer code.   One could even imagine that such a system could learn to accurately predict what sorts of definitions, conjectures, and theorems would be judged to be interesting by other mathematicians.  

We say this is ``troubling'' because it could generate the sort of input and output that we described above as essential for preserving the role of human mathematicians (or something indistinguishable from it), without humans having any role to play.  In such cases, one might think that there is no loss in full automation, since the automated system would be generating ``interesting'' results in a form that is comprehensible to humans.  But we think that even in this case, something crucial would be missing---though saying just what it is is more subtle.  We note two kinds of problem (and \citet{DeDeo} explores other, similar concerns).  The first is that the sort of scenario we are imagining, if truly devoid of human feedback, would involve a kind of problematic extrapolation.  Currently, mathematics evolves dynamically in response to new ideas, techniques, applications, and results.  That is, what mathematicians today would call interesting is not what mathematicians fifty or a hundred years ago would call interesting---and if the field is healthy, fifty years from now, mathematicians'  interests will have evolved further still.  LLMs, however, work by predicting next words based on frequencies in an existing body of text (including, perhaps, updates).  This can simulate many aspects of human discourse. 
 But making probabilistic judgments concerning what some fixed group of mathematicians would judge as interesting does not reflect the dynamics that has historically guided mathematics.  It seems to us something is lost.

Of course, one could also imagine a situation in which humans provide ongoing feedback to the AI mathematicians about what is interesting \emph{now}, in light of changing background conditions.  It might even be that mathematicians, responding to results generated by AI, would come to find new topics interesting, leading to new work by the AI, leading to new results, still---in each case, interesting within the evolving context.  But we would argue that even in this case, something would be missing.  This leads to the second problem, which is that the reasons the theorems are conjectured and proved would be the wrong reasons, and therefore they would not play the right role in the social activity of mathematics to count as contributions to that activity.  The picture, here, is one of human discovery of structures and relationships, leading to new ideas about how to investigate those structures and what significance they have for other things of prior interest, leading to yet more ideas about those structures and new ones.  The new ideas proceed from ever-deeper understanding of the old ones. This is an essentially creative activity, and a social one that stems from engagement with other mathematicians' work.  The process by which it occurs is ultimately what contributes most to our understanding of the mathematical subject matter, and, by extension, to those parts of the world that we seek to describe with it.

One can entertain more speculative possibilities, which perhaps could play more constructive roles.  For instance, consider a customized AI system, trained on both the body of mathematical knowledge and a particular mathematicians' style, previous work, ways of thinking, and preferred questions.  That system, engaged in a dynamic conversation, could play the role of a close collaborator or graduate student, proposing ideas, translating others' work, and drawing connections that the mathematician can then respond to.  Or consider a situation where an AI system develops a body of new theorems, but then also generates pedagogical materials -- textbooks, video lectures, question and answer sessions, and the like -- to teach the new work to humans, much as a human mathematician would.\footnote{Thank you to Tim Gowers for pushing us on this sort of possibility.}  It is easy to see how in this sort of scenario, an AI assistant could be helpful, just as a good teacher or student could be helpful.  Many mathematicians engage with others' work, at least at some points in their career, through proxies of just the sort that an AI system could presumably simulate.

Again, though, we do not dispute that AI or other automated systems can be useful to mathematicians.  Of course they can be.  Our point is that no degree of simulation of the processes of mathematics via generation of the text and symbols that mathematicians produce to communicate with one another about mathematical ideas can play the cognitive role that reasoning about mathematics can play.  This is true even if the simulated material is indistinguishable from the sorts of text that mathematicians themselves would produce.  At very best, what one can hope for would be text that generates, in the mind of a mathematician who reads it, new interpretive and creative processes analogous to those that reading work by another mathematician would produce---thereby stimulating new ideas.\footnote{An anonymous reviewer asks: surely a well-written AI-produced proof of a difficult conjecture would be much more valuable to mathematics than a lone genius who solves the same problem and never shares it, and so it cannot merely be what happens in some mathematician or other's mind that matters.  And yes, we agree: mathematics is a group activity.  Our point is that the AI-generated proof would be of no value at all, even if it is a perfectly sound proof, except insofar as it contributes to the group activity---and doing that happens only insofar as it stimulates the same kinds of generative engagement with human mathematicians that other human mathematicians do.}  But before one declares that this is enough, note that again the contribution to mathematics is what the mathematician produces when attempting to engage with the AI generated text, and not the text itself, which is essentially sterile.

\section*{Acknowledgments}
The second author was supported in part by NSF Grant No. DMS-1944862. We are grateful to Jeremy Avigad, Josh Batson, Silvia De Toffoli, Kenny Easwaran, Tim Gowers, Michael Harris, Harvey Lederman, Penelope Maddy, Alexa McLain, Colin McLarty, Toby Meadows, and Akshay Venkatesh for comments on a previous draft.

\bibliographystyle{elsarticle-harv}
\bibliography{proof} 

\end{document}